\documentclass[12pt,a4paper]{article}
\usepackage[latin 2]{inputenc}
\usepackage{amssymb}
\usepackage{latexsym}
\begin{document}

\begin{center}
{\large\textbf{On the theorem of Davenport and generalized Dedekind sums}}

\vspace{5mm}

\textbf{Bence Borda}

{\footnotesize Department of Mathematics, Rutgers University

110 Frelinghuysen Road, Piscataway, NJ-08854, USA

Email: \texttt{bordabence85@gmail.com}}

\vspace{5mm}

{\footnotesize \textbf{Keywords:} discrepancy, generalized Dedekind sum}

{\footnotesize \textbf{Mathematics Subject Classification (2010):} 11K38, 11F20}
\end{center}

\vspace{5mm}

\begin{abstract}

A symmetrized lattice of $2n$ points in terms of an irrational real number $\alpha$ is considered in the unit square, as in the theorem of Davenport. If $\alpha$ is a quadratic irrational, the square of the $L^2$ discrepancy is found to be $c( \alpha ) \log n + O \left( \log \log n \right)$ for a computable positive constant $c( \alpha )$. For the golden ratio $\varphi$, the value $\sqrt{c (\varphi ) \log n}$ yields the smallest $L^2$ discrepancy of any sequence of explicitly constructed finite point sets in the unit square. If the partial quotients $a_k$ of $\alpha$ grow at most polynomially fast, the $L^2$ discrepancy is found in terms of $a_k$ up to an explicitly bounded error term. It is also shown that certain generalized Dedekind sums can be approximated using the same methods. For a special generalized Dedekind sum with arguments $a, b$ an asymptotic formula in terms of the partial quotients of $\frac{a}{b}$ is proved.

\end{abstract}

\vspace{5mm}

\noindent\textbf{1. Introduction.} Consider an arbitrary finite set $A \subset [0,1]^2$ in the unit square. For any $x,y \in [0,1]$ let

\begin{equation}
S_A(x,y) = \left| A \cap \left( [0,x) \times [0,y) \right) \right|
\end{equation}

\noindent denote the number of elements of $A$ in the rectangle $[0,x) \times [0,y)$. A classical result of K. Roth [8] in the theory of discrepancy is that for any finite set $A \subset [0,1]^2$ we have

\begin{equation}
\int_0^1 \!\!\! \int_0^1 \left( S_A (x,y) - |A| xy \right)^2 \, \mathrm{d}x \, \mathrm{d}y > C \cdot \log |A|
\end{equation}

\noindent for some universal constant $C>0$. The square root of the left hand side of (2) is sometimes called the mean square discrepancy, or the $L^2$ discrepancy of the set $A$.

Several constructions for the set $A$ show that (2) is best possible up to a constant factor, the first of which is due to H. Davenport. For a positive integer $n$ and an irrational real number $\alpha$ consider the set

\begin{equation}
A = A( n , \alpha ) = \left\{ \left( \left\{ \pm k \alpha \right\} , \frac{k}{n} \right) : 1 \le k \le n \right\}
\end{equation}

\noindent of $2n$ points, where $\left\{ x \right\}$ denotes the fractional part of $x$. Davenport [5] shows that if $\alpha$ is badly approximable, i.e. $\inf_{m>0} m \left\| m \alpha \right\| >0$, where $\left\| x \right\|$ denotes the distance of $x$ from the nearest integer, then for $A$ as in (3) we have

\begin{equation}
\int_0^1 \!\!\! \int_0^1 \left( S_A (x,y) - |A| xy \right)^2 \, \mathrm{d}x \, \mathrm{d}y < C(\alpha ) \cdot \log |A|
\end{equation}

\noindent for some positive constant $C(\alpha )$ depending only on $\alpha$.

The purpose of this paper is to find the precise order of magnitude of the left hand side of (4), where $A$ is as in (3). We will work with a weaker assumption than $\alpha$ being badly approximable, however: we will assume that the continued fraction representation $\alpha = [a_0; a_1 , a_2 , \dots ]$ satisfies $a_k = O( k^d )$ for some constant $d \ge 0$. Note that $\alpha$ is badly approximable if and only if this condition holds with $d=0$. The motivation for this generality comes from the fact that the partial quotients of Euler's number $e$ satisfy $a_k = O (k)$. In fact, there is a class of transcendental numbers related to Euler's number $e$, including e.g. $e^{\frac{2}{n}}$ for every positive integer $n$, the partial quotients of which satisfy the same condition. Since there are very few classes of irrational numbers the continued fraction representations of which are explicitly known, we wanted our results to hold for as many of them as possible.

\vspace{5mm}

\noindent\textbf{2. The theorem of Davenport.} The original proof of Davenport [5] of (4) heavily uses the properties of the sequence $\left\| m \alpha \right\|$. The first step toward finding the precise order of magnitude of the left hand side of (4) is to isolate its dependence on $\left\| m \alpha \right\|$ as follows.

\vspace{5mm}

\noindent\textsc{Theorem 1.} \textit{Let $\alpha$ be an irrational real number. Suppose its continued fraction representation $\alpha = [a_0; a_1 , a_2 , \dots ]$ satisfies $a_k = O(k^d)$ for some constant $d \ge 0$. For $A$ as in (3) we have}

\[ \int_0^1 \!\!\! \int_0^1 \left( S_A (x,y) - |A| xy \right)^2 \, \mathrm{d}x \, \mathrm{d}y = \sum_{m=1}^n \frac{1}{4 \pi^4 m^2 \left\| m \alpha \right\|^2} + O \left( \log^{2d} n \log \log n \right) , \]

\noindent\textit{as $n \to \infty$. The implied constant depends only on $\alpha$.}

\vspace{5mm}

The proof of Theorem 1 is based on the observation that $S_A(x,y)$ is constant on the horizontal stripes $\frac{k}{n} < y \le \frac{k+1}{n}$. This enables us to first integrate $\left( S_A (x,y) - |A|xy \right)^2$ with respect to $y$ on $[0,1]$, then apply the Parseval formula on the function obtained of the variable $x$. The same idea is used in the original proof of Davenport [5]. To isolate the specific Diophantine sum in Theorem 1, we will need delicate estimates of certain trigonometric sums, however, which are not present in the original paper.

\vspace{5mm}

\noindent\textbf{3. Generalized Dedekind sums.} A Diophantine sum similar to the one in Theorem 1 appears in the study of generalized Dedekind sums. Following the notation of [6] let $B_n (x)$ denote the $n$th Bernoulli polynomial, defined e.g. recursively as $B_0(x)=1$,

\begin{equation}
B_n' (x) = nB_{n-1}(x), \qquad \int_0^1 B_n(x) \, \mathrm{d}x=0
\end{equation}

\noindent for every positive integer $n$. For relatively prime positive integers $a,b$ and positive integers $p,q$ we define the generalized (inhomogeneous) Dedekind sum $s_{p,q}(a,b)$ as

\begin{equation}
s_{p,q} (a,b) = \sum_{k=1}^{b-1} B_p \left( \frac{k}{b} \right) B_{q} \left( \left\{ \frac{ak}{b} \right\} \right) .
\end{equation}

In the special case, when $p=q$ is even, the generalized Dedekind sum can be approximated by a Diophantine sum similar to the one in Theorem 1, with a \textit{rational} number $\alpha=\frac{a}{b}$.

\vspace{5mm}

\noindent\textsc{Theorem 2.} \textit{Let $a,b$ be relatively prime positive integers, and let $p \ge 2$ be an even integer. Then}

\[ s_{p,p} (a,b) = \frac{2 \left( p! \right)^2}{(2 \pi)^{2p} b^{p-1}} \left( \sum_{m=1}^{b-1} \frac{1}{m^p \left\| m \frac{a}{b} \right\|^p} + E \right) , \]

\noindent\textit{where the error satisfies $0<E<5 \cdot 2^p$.}

\vspace{5mm}

Applying integration by parts $p$ times, and using the recursive definition (5) of the Bernoulli polynomials it is easy to see that for any integer $m \neq 0$ their Fourier coefficients are

\[ \int_0^1 B_p (x) e^{-2 \pi i m x} \, \mathrm{d}x = \frac{- p!}{(2 \pi i m )^p} . \]

\noindent Thus the Fourier series of $B_p \left( \left\{ x \right\} \right)$ is particularly simple:

\begin{equation}
B_p \left( \left\{ x \right\} \right) = \sum_{m \neq 0} \frac{-p!}{(2 \pi i m )^p} e^{2 \pi i m x}
\end{equation}

\noindent for any real number $x$ and any integer $p \ge 2$. Moreover, the series is absolutely convergent. The main idea in the proof of Theorem 2 is to replace $B_p$ in the definition (6) of $s_{p,p} (a,b)$ by its Fourier series (7).

\newpage

\noindent\textbf{4. The Diophantine sum.} We now turn our attention to the Diophantine sum

\begin{equation}
\sum_{m=1}^n \frac{1}{m^p \left\| m \alpha \right\|^p} .
\end{equation}

\noindent J. Beck [2] studies the similar sum

\[ \sum_{m=1}^n \frac{1}{m^2 \sin^2 (m \pi \alpha )} \]

\noindent in the special case, when $\alpha$ is a quadratic irrational. Since

\[ \sum_{m=1}^n \frac{1}{m^2 \sin^2 (m \pi \alpha )} = \sum_{m=1}^n \frac{1}{\pi^2 m^2 \left\| m \alpha \right\|^2} + O(1) \]

\noindent with an absolute implied constant, his results are directly applicable. It is shown ([2] Proposition 3.2) that if $\alpha$ is a quadratic irrational, then

\begin{equation}
\sum_{m=1}^n \frac{1}{4 \pi^4 m^2 \left\| m \alpha \right\|^2} = c (\alpha ) \log n + O(1)
\end{equation}

\noindent for some positive constant $c(\alpha )$ depending only on $\alpha$, where the implied constant depends only on $\alpha$ as well. Several methods are known ([2] Chapter 3) to compute $c( \alpha )$, e.g. ([2] 3.92 , 3.102) we have

\[ c(\sqrt{2}) = \frac{1}{24 \sqrt{2} \log \left( 1 + \sqrt{2} \right)} , \qquad c(\sqrt{3}) = \frac{1}{12 \sqrt{3} \log \left( 2 + \sqrt{3} \right)} , \]

\[ c \left( \frac{1+\sqrt{5}}{2} \right) = \frac{1}{30 \sqrt{5} \log \left( \frac{1+ \sqrt{5}}{2} \right)} . \]

\noindent A rather general result ([2] 3.2.1) is that for any prime $p \equiv 3 \pmod{4}$ such that the class number $h(p)=1$ we have

\[ c (\sqrt{p} ) = \frac{p}{\pi^4 \log \eta_p} \zeta_K (2) , \]

\noindent where $\zeta_K$ is the Dedekind zeta function, while $\eta_p$ is the fundamental unit of the real quadratic field $K=\mathbb{Q} \left( \sqrt{p} \right)$.

We now offer a way to estimate the Diophantine sum (8) in terms of the continued fraction representation of $\alpha$.

\vspace{5mm}

\noindent\textsc{Theorem 3.} \textit{Let $\alpha$ be an arbitrary real number, and consider its (finite or infinite) continued fraction representation $\alpha = [a_0; a_1 , a_2 , \dots ]$. Let $\frac{p_k}{q_k}=[a_0 ; a_1 , \dots , a_{k-1}]$ denote the convergents to $\alpha$. For any real number $p>1$ and any positive integer $\ell$ (which is at most the number of partial quotients in the case of a rational $\alpha$) we have}

\[ \left( \sum_{0<m<q_{\ell}} \frac{1}{m^p \left\| m \alpha \right\|^p} \right)^{\frac{1}{p}} = \left( \sum_{0<k< \ell} \zeta (2p) a_k^p \right)^{\frac{1}{p}} + O \left( \ell^{\frac{1}{p}} \right) , \]

\noindent\textit{where $\zeta$ is the Riemann zeta function. The implied constant is at most $6^p \frac{4p^2}{(p-1)^2}$.}

\vspace{5mm}

In the proof of Theorem 3 we will decompose the Diophantine sum into two parts. For the terms where $m$ is an integral multiple $m=a q_k$ of a convergent denominator $q_k$ for some $0<k< \ell$ and $a>0$ we have

\[ \frac{1}{(a q_k)^p \left\| a q_k \right\|^p} \approx \frac{a_k^p}{a^{2p}} . \]

\noindent We will obtain the main term by summing over every positive integer $a$ and over $0<k< \ell$. The terms where $m$ is not of the form $m=a q_k$ will be treated as an error term.

Note that if the partial quotients $a_k$ of $\alpha$ are bounded, in particular if $\alpha$ is a quadratic irrational, then the main term and the error term in Theorem 3 have the same order of magnitude, making the result useless. If, however,

\[ \frac{1}{\ell} \sum_{0<k< \ell} a_k^p \to \infty , \]

\noindent in particular if $\alpha$ is Euler's number $e$, then Theorem 3 evaluates the Diophantine sum (8) up to an error of smaller order of magnitude.

\vspace{5mm}

\noindent\textbf{5. Conclusions.} We can easily combine Theorem 1, and the result (9) of J. Beck or Theorem 3 with $p=2$ to obtain the following.

\vspace{5mm}

\noindent\textsc{Corollary 4.} \textit{Let $\alpha$ be an irrational real number, and let $A$ be as in (3).}

\begin{enumerate}
\item[(i)] \textit{If $\alpha$ is a quadratic irrational, then}

\[ \int_0^1 \!\!\! \int_0^1 \left( S_A (x,y) - |A| xy \right)^2 \, \mathrm{d}x \, \mathrm{d}y = c(\alpha ) \log n + O \left( \log \log n \right) \]

\textit{for some positive constant $c (\alpha )$, where $c(\alpha )$ and the implied constant depend only on $\alpha$.}

\item[(ii)] \textit{Suppose the continued fraction representation $\alpha = [a_0 ; a_1 , a_2 , \dots ]$ satisfies $a_k = O\left( k^d \right)$ for some constant $d \ge 0$. Let $\frac{p_k}{q_k} = [a_0 ; a_1 , a_2 , \dots , a_{k-1}]$ denote the convergents to $\alpha$. If $q_{\ell} \le n < q_{\ell +1}$ then}

\[ \int_0^1 \!\!\! \int_0^1 \left( S_A (x,y) - |A| xy \right)^2 \, \mathrm{d}x \, \mathrm{d}y = \]

\[ \frac{1}{360} \sum_{k=1}^{\ell} a_k^2 + O \left( \log^{d+1} n + \log^{2d} n \log \log n \right) . \]

\textit{The implied constant depends only on $\alpha$.}
\end{enumerate}

\begin{flushright}
$\square$
\end{flushright}

\vspace{5mm}

Intuitively one expects that the smaller the partial quotients of $\alpha$ are, the smaller the left hand side of (4) is. The most extreme case is that of the golden ratio $\alpha = \frac{1+\sqrt{5}}{2} = [1;1,1, \dots ]$, for which we have ([2] 3.102)

\[ c \left( \frac{1+\sqrt{5}}{2} \right) = \frac{1}{30 \sqrt{5} \log \left( \frac{1+ \sqrt{5}}{2} \right)} \approx 0.030978 . \]

\noindent Note that ([3] 2.61) numerical computations have already yielded the constant 0.030978 for an essentially identical construction for the set $A$, although it has neither been supported by a rigorous proof, nor has it been identified as an explicit expression. According to [7], as of 2016 every other known construction for a sequence of finite sets $A \subset [0,1]^2$ satisfies

\[ \liminf_{|A| \to \infty} \frac{1}{\log |A|} \int_0^1 \!\!\! \int_0^1 \left( S_A (x,y) - |A| xy \right)^2 \, \mathrm{d}x \, \mathrm{d}y > 0.032. \]

Now let $a,b$ be relatively prime positive integers, let $p \ge 2$ be an even integer, and consider the generalized Dedekind sum $s_{p,p} (a,b)$. From a computational point of view, Theorem 2 is not an effective way of approximating $s_{p,p} (a,b)$: both definition (6) and the formula in Theorem 2 require the computation of a sum of $b-1$ terms. If we consider the continued fraction representation $\frac{a}{b} = [a_0 ; a_1 , \dots , a_{\ell} ]$, however, we have that the last convergent is $\frac{p_{\ell +1}}{q_{\ell +1}} = \frac{a}{b}$, therefore we can combine Theorem 2 and Theorem 3 with $q_{\ell +1}=b$ to obtain

\begin{equation}
s_{p,p} (a,b) = \frac{2 (p!)^2 \zeta (2p)}{(2 \pi )^{2p} b^{p-1}} \sum_{k=1}^{\ell} a_k^p \left( 1 + O \left( \left( \frac{1}{\ell} \sum_{k=1}^{\ell} a_k^p \right)^{- \frac{1}{p}} \right) \right)
\end{equation}

\noindent with an implied constant depending only on $p$. Note that (10) requires the computation of a sum of only $\ell = O \left( \log b \right)$ terms. On the other hand, (10) is only useful if the power mean of the partial quotients satisfies

\[ \left( \frac{1}{\ell} \sum_{k=1}^{\ell} a_k^p \right)^{\frac{1}{p}} \to \infty . \]

It is interesting to compare (10) to Barkan's evaluation [1] of $s_{1,1} (a,b)$, which roughly states

\begin{equation}
s_{1,1} (a,b) = \frac{1}{12} \sum_{k=1}^{\ell} (-1)^{k+1} a_k + O(1) .
\end{equation}

\noindent The proof of (11) is based on a reciprocity formula for $s_{1,1} (a,b)$, and in fact gives the exact value without an error term. While various reciprocity formulas are known for the generalized Dedekind sum $s_{p,q} (a,b)$ as well (see e.g. [6]), these formulas do not yield an evaluation similar to (11). Although asymptotic formulas for $s_{p,q} (a,b)$ are known (e.g. [6]), (10) seems to be the first asymptotic formula in terms of the partial quotients of $\frac{a}{b}$.

\vspace{5mm}

\noindent\textbf{6. The proofs of the theorems.} We start by stating the facts about continued fractions to be used in the proofs. We follow the conventions of Cassels [4]. The finite or infinite continued fraction representation of an arbitrary real number $\alpha$ will be denoted by $\alpha = [a_0 ; a_1 , a_2 , \dots ]$. In the case of a rational $\alpha$ it will not matter which of the two possible representations is chosen. The convergents to $\alpha$ will be denoted by

\[ \frac{p_k}{q_k} = [a_0 ; a_1 , \dots , a_{k-1}] = a_0 + \frac{1}{a_1 + \frac{1}{a_2 + \frac{1}{\cdots + \frac{1}{a_{k-1}}}}} . \]

\noindent The properties of continued fractions to be used are the following.

\vspace{5mm}

\noindent\textsc{Proposition 5.} \textit{For an arbitrary real number $\alpha$ we have:}

\begin{enumerate}
\item[(i)] \textit{If either $k \ge 2$, or $k=1$ and $a_1>1$, then $\frac{1}{q_{k+1}+q_k} \le \left\| q_k \alpha \right\| \le \frac{1}{q_{k+1}}$}.
\item[(ii)] \textit{We have $q_1=1$, $q_2=a_1$ and $q_{k+1}=a_k q_k + q_{k-1}$ for any $k \ge 2$.}
\item[(iii)] \textit{For any $k \ge 2$ we have $p_k q_{k-1} - q_k p_{k-1} = (-1)^k$, and the numbers $p_k, q_k$ are relatively prime.}
\item[(iv)] \textit{For any $k \ge 1$ we have $\mathrm{sign} \left( q_k \alpha - p_k \right) = (-1)^{k+1}$.}
\item[(v)] \textit{For any $k \ge 1$ we have $q_1 + q_2 + \cdots + q_k \le 3 q_k$.}
\item[(vi)] \textit{For any $k \ge 2$ we have}

\[\sum_{0<m<q_k} \frac{1}{\left\| m \alpha \right\|} \le 8 q_k \log_2 (2 q_k) . \]

\item[(vii)] \textit{For any $k \ge 2$ and $p >1$ we have}

\[ \sum_{0<m<q_k} \frac{1}{\left\| m \alpha \right\|^p} \le \frac{4^{p+1}}{2^p -2} q_k^p . \]

\end{enumerate}

\vspace{5mm}

\noindent\textbf{Proof:} (i)-(iv) See e.g. [4] Chapter I.

\vspace{5mm}

\noindent (v) Summing the recurrence in (ii) we get

\[ q_k + q_{k-1}-q_2-q_1 = \sum_{\ell =2}^{k-1} a_{\ell} q_{\ell} \ge \sum_{\ell =2}^{k-1} q_{\ell} , \]

\[ 2 q_k \ge q_{k} + q_{k-1} - q_2 \ge \sum_{\ell =1}^{k-1} q_{\ell} . \]

\noindent Adding $q_k$ finishes the proof.

\vspace{5mm}

\noindent (vi)-(vii) For any $0<m<q_k$ we have $\left\| m \alpha \right\| \ge \left\| q_{k-1} \alpha \right\| > \frac{1}{2 q_k}$. For any integer $n \ge 0$ consider the set

\[ A_n = \left\{ 0 < m < q_k : 2^n \frac{1}{2 q_k} \le \left\| m \alpha \right\| < 2^{n+1} \frac{1}{2 q_k} \right\} . \]

\noindent If $2^n \frac{1}{2 q_k} > \frac{1}{2}$, i.e. if $n > \log_2 q_k$, then $A_n = \emptyset$.

For every $m \in A_n$ consider the point in $\left[ - \frac{1}{2} , \frac{1}{2} \right)$ equivalent to $m \alpha$ modulo 1. All these points belong to the interval $\left( - 2^{n+1} \frac{1}{2 q_k} , 2^{n+1} \frac{1}{2 q_k} \right)$ of length $\frac{2^{n+1}}{q_k}$. On the other hand, the distance of any two points is larger, than $\frac{1}{2 q_k}$. Indeed, otherwise there would exist integers $0<m<m'<q_k$ such that $\left\| (m'-m) \alpha \right\| \le \frac{1}{2 q_k}$. Hence by the pigeonhole principle we have $|A_n| \le 2^{n+2}$.

First consider

\[ \sum_{0<m<q_k} \frac{1}{\left\| m \alpha \right\|} = \sum_{0 \le n \le \log_2 q_k} \sum_{m \in A_n} \frac{1}{\left\| m \alpha \right\|} \le \sum_{0 \le n \le \log_2 q_k} \frac{2q_k}{2^n} 2^{n+2} \le 8 q_k \log_2 (2 q_k ). \]

\noindent Finally,

\[ \sum_{0<m<q_k} \frac{1}{\left\| m \alpha \right\|^p} = \sum_{0 \le n \le \log_2 q_k} \sum_{m \in A_n} \frac{1}{\left\| m \alpha \right\|^p} \le \sum_{n=0}^{\infty} \frac{(2q_k)^p}{2^{pn}} 2^{n+2} = \frac{4^{p+1}}{2^p -2} q_k^p . \]

\begin{flushright}
$\square$
\end{flushright}

\vspace{5mm}

\noindent\textbf{Proof of Theorem 1:} Elementary calculation shows that for any real number $S$ and any integer $0 \le k \le n-1$ we have

\begin{equation}
\int_{\frac{k}{n}}^{\frac{k+1}{n}} \left( S - 2nxy \right)^2  \, \mathrm{d}y = \frac{1}{n} \left( S-2kx \right)^2 - \frac{2}{n} \left( S-2kx \right) x + \frac{4}{3n} x^2 .
\end{equation}

\noindent From the definition (1) of $S_A (x,y)$ it is clear that for any $0 \le x \le 1$ and any $\frac{k}{n} < y \le \frac{k+1}{n}$ we have $S_A (x,y) = S_A \left( x , \frac{k+1}{n} \right)$. To compute the double integral in the theorem, we can thus substitute $S=S_A \left( x , \frac{k+1}{n} \right)$ in (12), integrate with respect to $x$ on $[0,1]$, and sum over $0 \le k \le n-1$. Let us introduce a notation for the first two terms obtained:

\[ M = \frac{1}{n} \sum_{k=0}^{n-1} \int_0^1 \left( S_A \left( x , \frac{k+1}{n} \right) - 2kx \right)^2 \, \mathrm{d}x , \]

\[ L = - \frac{2}{n} \sum_{k=0}^{n-1} \int_0^1 \left( S_A \left( x , \frac{k+1}{n} \right) - 2kx \right) x \, \mathrm{d}x . \]

\noindent Since the contribution of $\frac{4}{3n} x^2$ in (12) is clearly bounded, we have

\begin{equation}
\int_0^1 \!\!\! \int_0^1 \left( S_A (x,y) - |A| xy \right)^2 \, \mathrm{d}x \, \mathrm{d}y = M + L + O \left( 1 \right) .
\end{equation}

Now we compute the main term $M$ and estimate $L$. From the definition (1) of $S_A$ we can see that

\[ S_A \left( x , \frac{k+1}{n} \right) - 2kx = \sum_{j=1}^k \bigg( \chi_{[0,x)} \left( \left\{ j \alpha \right\} \right) + \chi_{[0,x)} \left( \left\{ - j \alpha \right\} \right) -2x \bigg) , \]

\noindent where $\chi$ denotes the characteristic function of a set. Elementary integration shows

\[ \int_0^1 \bigg( \chi_{[0,x)} \left( \left\{ j \alpha \right\} \right) + \chi_{[0,x)} \left( \left\{ - j \alpha \right\} \right) -2x \bigg) \, \mathrm{d}x =0 , \]

\[ \int_0^1 \bigg( \chi_{[0,x)} \left( \left\{ j \alpha \right\} \right) + \chi_{[0,x)} \left( \left\{ - j \alpha \right\} \right) -2x \bigg) e^{-2 \pi i m x} \, \mathrm{d}x = \frac{\cos (2 m j \pi \alpha)}{\pi i m} \]

\noindent for any integer $m \neq 0$. Summing over $1 \le j \le k$ we obtain the Fourier coefficients

\[ \int_{0}^1 \left( S_A \left( x , \frac{k+1}{n} \right) - 2kx \right) \, \mathrm{d}x =0 , \]

\[ \int_0^1 \left( S_A \left( x , \frac{k+1}{n} \right) - 2kx \right) e^{- 2 \pi i m x} \, \mathrm{d}x = \frac{\sin \left( (2k+1) m \pi \alpha \right)}{2 \pi i m \sin (m \pi \alpha)} + O \left( \frac{1}{|m|} \right) \]

\noindent for any integer $m \neq 0$. Note that $|\sin (m \pi \alpha )| \ge 2 \left\| m \alpha \right\|$.

Since the Fourier coefficients of $x$ are $O \left( \frac{1}{|m|} \right)$, the Parseval formula yields

\begin{equation}
L = O \left( \sum_{m=1}^{\infty} \frac{1}{m^2 \left\| m \alpha \right\|} \right) .
\end{equation}

\noindent Applying the Parseval formula again we obtain

\begin{equation}
M = \frac{1}{n} \sum_{k=0}^{n-1} \sum_{m=1}^{\infty} \frac{\sin^2 \left( (2k+1) m \pi \alpha \right)}{2 \pi^2 m^2 \sin^2 (m \pi \alpha )} +  O \left( \sum_{m=1}^{\infty} \frac{1}{m^2 \left\| m \alpha \right\|} \right) .
\end{equation}

\noindent The error terms in (14) and (15) can be estimated as

\begin{equation}
\sum_{m=1}^{\infty} \frac{1}{m^2 \left\| m \alpha \right\|} = \sum_{k=1}^{\infty} \sum_{q_k \le m < q_{k+1}} \frac{1}{m^2 \left\| m \alpha \right\|} \le \sum_{k=1}^{\infty} \frac{1}{q_k^2} \sum_{0 < m < q_{k+1}} \frac{1}{\left\| m \alpha \right\|} .
\end{equation}

\noindent Using Proposition 5 (vi), $q_{k+1} = O \left( a_k q_k \right) = O \left( k^d q_k \right)$, and the fact that $q_k$ is at least as big as the $k$th Fibonacci number, it is easy to see that (16) is $O(1)$. Thus (13)-(15) yield

\begin{equation}
\int_0^1 \!\!\! \int_0^1 \left( S_A (x,y) - |A| xy \right)^2 \, \mathrm{d}x \, \mathrm{d}y = \frac{1}{n} \sum_{k=0}^{n-1} \sum_{m=1}^{\infty} \frac{\sin^2 \left( (2k+1) m \pi \alpha \right)}{2 \pi^2 m^2 \sin^2 (m \pi \alpha )} + O(1).
\end{equation}

We now estimate the tail of the infinite series in (17). Note that

\[ \frac{\sin^2 \left( (2k+1) m \pi \alpha \right)}{\sin^2 (m \pi \alpha )} = \left| \sum_{j=0}^{2k} e^{2 \pi i j m \alpha} \right|^2 = 2k+1 + \sum_{{{0 \le j_1, j_2 \le 2k}\atop{j_1 \neq j_2}}} e^{2 \pi i (j_1 - j_2) m \alpha} . \]

\noindent Hence the partial sums satisfy

\[ c_{\ell} = \sum_{m=1}^{\ell} \frac{\sin^2 \left( (2k+1) m \pi \alpha \right)}{\sin^2 (m \pi \alpha )} = (2k+1) \ell + \sum_{{{0 \le j_1, j_2 \le 2k}\atop{j_1 \neq j_2}}} \frac{e^{2 \pi i (j_1 - j_2) \ell \alpha}-1}{1-e^{-2 \pi i (j_1 - j_2 ) \alpha }} = \]

\[ O \left( k \ell + \sum_{{{0 \le j_1, j_2 \le 2k}\atop{j_1 \neq j_2}}} \frac{1}{\left\| (j_1 - j_2) \alpha \right\|} \right) = O \left( k \ell + k \sum_{j=1}^{2k} \frac{1}{\left\| j \alpha \right\|} \right) . \]

\noindent Proposition 5 (vi) and $a_k=O (k^d)$ thus yield $c_{\ell} = O \left( k \ell + k^2 \log^{d+1} k \right)$. Applying summation by parts on the infinite series in (17) starting at $m = \lfloor k \sqrt{\log k} \rfloor$ we get

\[ \sum_{m= \lfloor k \sqrt{\log k} \rfloor}^{\infty} \frac{1}{2 \pi^2 m^2} \cdot \frac{\sin^2 \left( (2k+1) m \pi \alpha \right)}{\sin^2 (m \pi \alpha )} = \]

\[ - c_{\lfloor k \sqrt{\log k} \rfloor -1} \cdot \frac{1}{2 \pi^2 \lfloor k \sqrt{\log k} \rfloor^2} + \sum_{m= \lfloor k \sqrt{\log k} \rfloor}^{\infty} c_m \cdot \left( \frac{1}{2 \pi^2 m^2} - \frac{1}{2 \pi^2 (m+1)^2} \right) = \]

\[ O \left( \log^d k + \sum_{m=\lfloor k \sqrt{\log k} \rfloor}^{\infty} \frac{k m + k^2 \log^{d+1} k}{m^3} \right) = O \left( \log^d k \right) . \]

\noindent Therefore we can replace the infinite series in (17) by a finite sum to obtain

\[ \int_0^1 \!\!\! \int_0^1 \left( S_A (x,y) - |A| xy \right)^2 \, \mathrm{d}x \, \mathrm{d}y = \]

\begin{equation}
\frac{1}{n} \sum_{k=0}^{n-1} \sum_{1 \le m \le n \sqrt{\log n}} \frac{\sin^2 \left( (2k+1) m \pi \alpha \right)}{2 \pi^2 m^2 \sin^2 (m \pi \alpha )} + O( \log^d n ).
\end{equation}

Let us switch the order of summation in (18), and use the trigonometric identity

\[ \sum_{k=0}^{n-1} \sin^2 \left( (2k+1) x \right) = \frac{n}{2} - \frac{\sin (4 n x)}{4 \sin (2x)} \]

\noindent with $x= m \pi \alpha$ to get

\[ \int_0^1 \!\!\! \int_0^1 \left( S_A (x,y) - |A| xy \right)^2 \, \mathrm{d}x \, \mathrm{d}y = \sum_{1 \le m \le n \sqrt{\log n}} \frac{1}{4 \pi^2 m^2 \sin^2 (m \pi \alpha )} + \]

\begin{equation}
O \left( \frac{1}{n} \sum_{1 \le m \le n \sqrt{\log n}} \frac{\left| \sin (4 n m \pi \alpha ) \right|}{m^2 \left| \sin (2 m \pi \alpha ) \right|^3} + \log^d n \right) .
\end{equation}

We first estimate the terms $1 \le m \le \frac{n}{\log^{3d} n}$ in the error. We have

\[ \frac{1}{n} \sum_{1 \le m \le \frac{n}{\log^{3d}n}} \frac{\left| \sin (4 n m \pi \alpha ) \right|}{m^2 \left| \sin (2 m \pi \alpha ) \right|^3} = O \left( \sum_{1 \le m \le \frac{2n}{\log^{3d} n}} \frac{1}{n m^2 \left\| m \alpha \right\|^3} \right) = \]

\begin{equation}
O \left( \sum_{{{k}\atop{q_k \le \frac{2n}{\log^{3d}n}}}} \sum_{q_k \le m < q_{k+1}} \frac{1}{n m^2 \left\| m \alpha \right\|^3} \right) = O \left( \sum_{{{k}\atop{q_k \le \frac{2n}{\log^{3d}n}}}} \frac{1}{n q_k^2} \sum_{0 < m < q_{k+1}} \frac{1}{\left\| m \alpha \right\|^3} \right) .
\end{equation}

\noindent Using Proposition 5 (vii), $q_{k+1}^3 = O \left( \log^{3d} n \cdot q_k^3 \right)$ and Proposition 5 (v) one can see that (20) is $O(1)$. We can estimate the terms $\frac{n}{\log^{3d}n} \le m \le n \sqrt{\log n}$ in the error of (19) by applying $|\sin (4 n m \pi \alpha)| \le 2n \left|\sin (2 m \pi \alpha) \right|$ to get

\[ \int_0^1 \!\!\! \int_0^1 \left( S_A (x,y) - |A| xy \right)^2 \, \mathrm{d}x \, \mathrm{d}y = \sum_{m=1}^n \frac{1}{4 \pi^2 m^2 \sin^2 (m \pi \alpha )} + \]

\begin{equation}
O \left( \sum_{\frac{n}{\log^{3d} n} \le m \le 2 n \sqrt{\log n}} \frac{1}{m^2 \left\| m \alpha \right\|^2} + \log^d n \right).
\end{equation}

Finally, we need to estimate the error in (21). We can use Proposition 5 (vii) again to estimate the terms as $m$ runs between two consecutive convergent denominators as

\[ \sum_{q_k \le m < q_{k+1}} \frac{1}{m^2 \left\| m \alpha \right\|^2} \le \frac{1}{q_k^2} \sum_{0<m<q_{k+1}} \frac{1}{\left\| m \alpha \right\|^2} = O \left( k^{2d} \right) = O \left( \log^{2d} n \right) . \]

\noindent The recurrence in Proposition 5 (ii) yields $\frac{q_{k+2}}{q_k} = \frac{a_{k+1}q_{k+1} + q_k}{q_k} \ge 2$, which in turn shows that the number of convergent denominators which fall in the interval $\left[ \frac{n}{\log^{3d}n} , 2n \sqrt{\log n} \right]$ is $O \left( \log \log n \right)$. Thus the error in (21) is $O \left( \log^{2d} n \log \log n \right)$, hence

\[ \int_0^1 \!\!\! \int_0^1 \left( S_A (x,y) - |A| xy \right)^2 \, \mathrm{d}x \, \mathrm{d}y = \sum_{m=1}^n \frac{1}{4 \pi^2 m^2 \sin^2 (m \pi \alpha )} + O \left( \log^{2d} n \log \log n \right) . \]

Finally, note that

\[ \sum_{m=1}^n \frac{1}{4 \pi^2 m^2 \sin^2 (m \pi \alpha )} = \sum_{m=1}^n \frac{1}{4 \pi^4 m^2 \left\| m \alpha \right\|^2} + O(1). \]

\begin{flushright}
$\square$
\end{flushright}

\vspace{5mm}

\noindent\textbf{Proof of Theorem 2:} Since $B_p (0)=0$, we can change the lower limit of summation in the definition (6) of $s_{p,p} (a,b)$ to $k=0$. Substituting the Fourier series (7) of $B_p \left( \left\{ x \right\} \right)$ in (6) we obtain

\[ s_{p,p} (a,b) = \sum_{m \neq 0} \sum_{\ell \neq 0} \frac{(p!)^2}{(2 \pi)^{2p} m^p \ell^p} \sum_{k=0}^{b-1} e^{2 \pi i \frac{\ell + a m}{b} k} . \]

\noindent Note that the inner sum is zero whenever $b \nmid \ell + a m$, therefore

\begin{equation}
s_{p,p} (a,b) = \frac{(p!)^2 b}{(2 \pi )^{2p}} \sum_{m \neq 0} \sum_{{{\ell \neq 0}\atop{b \mid \ell + am}}} \frac{1}{m^p \ell^p}.
\end{equation}

It is easy to see that the sum of all the terms in (22) such that $b \mid m$ is

\[ \frac{(p!)^2}{(2 \pi)^{2p} b^{2p-1}} 4 \zeta (p)^2 . \]

\noindent For every term in (22) such that $b \nmid m$ we have that the inner sum is over integers $\ell = jb -am$, as $j$ runs in $\mathbb{Z}$. Thus

\begin{equation}
s_{p,p} (a,b) = \frac{(p!)^2}{(2 \pi )^{2p} b^{p-1}} \left( \sum_{b \nmid m} \sum_{j \in \mathbb{Z}} \frac{1}{m^p \left( j- \frac{a}{b} m \right)^p} + \frac{4}{b^p} \zeta (p)^2 \right) .
\end{equation}

Note that for any $c \not\in \mathbb{Z}$ we have

\[ \frac{1}{\left\| c \right\|^p} \le \sum_{j \in \mathbb{Z}} \frac{1}{(j-c)^p} \le \frac{1}{\left\| c \right\|^p} + \left( \sum_{k=1}^{\infty} \frac{1}{k^p} + \sum_{k=0}^{\infty} \frac{1}{\left( k + \frac{1}{2} \right)^p} \right) = \frac{1}{\left\| c \right\|^p} + 2^p \zeta (p) . \]

\noindent Applying this with $c = \frac{a}{b}m$ in (23) we get

\[ s_{p,p} (a,b) = \frac{(p!)^2}{(2 \pi )^{2p} b^{p-1}} \left( \sum_{b \nmid m} \frac{1}{m^p \left\| m \frac{a}{b} \right\|^p} + A \right) \]

\noindent for some

\[ 0 < A < 2^{p+1} \zeta (p)^2 + \frac{4}{b^p} \zeta (p) . \]

The terms indexed by $m$ and $-m$ are equal. Thus

\[ s_{p,p} (a,b) = \frac{2(p!)^2}{(2 \pi )^{2p} b^{p-1}} \left( \sum_{m=1}^{b-1} \frac{1}{m^p \left\| m \frac{a}{b} \right\|^p} + \sum_{k=1}^{\infty} \sum_{kb < m < (k+1)b} \frac{1}{m^p \left\| m \frac{a}{b} \right\|^p} + \frac{A}{2} \right) . \]

\noindent Here we have

\[ \sum_{kb < m < (k+1)b} \frac{1}{m^p \left\| m \frac{a}{b} \right\|^p} \le \frac{1}{k^p b^p} \sum_{kb < m < (k+1)b} \frac{1}{\left\| \frac{ma}{b} \right\|^p} . \]

\noindent Since $(a,b)=1$, as $m$ runs in $kb < m < (k+1)b$, the numbers $ma$ fall into each nonzero residue class modulo $b$ exactly once. Therefore the sum is

\[ \frac{1}{k^p b^p} \sum_{j=1}^{b-1} \frac{1}{\left\| \frac{j}{b} \right\|^p} \le \frac{2}{k^p} \zeta (p) . \]

\noindent Hence

\[ s_{p,p} (a,b) = \frac{2(p!)^2}{(2 \pi )^{2p} b^{p-1}} \left( \sum_{m=1}^{b-1} \frac{1}{m^p \left\| m \frac{a}{b} \right\|^p} + E \right) \]

\noindent for some

\[ 0 < E < 2 \zeta (p)^2 + 2^p \zeta (p)^2 + \frac{2}{b^p} \zeta (p) < 5 \cdot 2^p . \]

\begin{flushright}
$\square$
\end{flushright}

\newpage

\noindent\textbf{Proof of Theorem 3:} First suppose $\ell \ge 2$ and $q_{\ell} \ge 2$, and consider the sum

\begin{equation}
\sum_{q_{\ell} \le m < q_{\ell +1}} \frac{1}{m^p \left\| m \alpha \right\|^p} .
\end{equation}

\noindent Let us introduce the notation $\varepsilon_{\ell} = q_{\ell} \alpha - p_{\ell}$. Then

\[ \left\| m \alpha \right\| = \left\| \frac{m p _{\ell}}{q_{\ell}} + \frac{m \varepsilon_{\ell}}{q_{\ell}} \right\| . \]

\noindent We will decompose the sum (24) using the index sets

\[ A= \left\{ q_{\ell} \le m < q_{\ell +1} : m p_{\ell} \not\equiv 0 , (-1)^{\ell} \pmod{q_{\ell}} \right\} , \]

\[ B= \left\{ q_{\ell} \le m < q_{\ell +1} : m p_{\ell} \equiv (-1)^{\ell} \pmod{q_{\ell}} \right\} , \]

\[ C= \left\{ q_{\ell} \le m < q_{\ell +1} : m p_{\ell} \equiv 0 \pmod{q_{\ell}} \right\} . \]

Consider the sum over $m \in A$ first. The assumption $q_{\ell} \ge 2$ and Proposition 5 (i) imply that for any $q_{\ell} \le m < q_{\ell +1}$

\[ \left| \frac{m \varepsilon_{\ell}}{q_{\ell}} \right| = \frac{m \left\| q_{\ell} \alpha \right\|}{q_{\ell}} < \frac{1}{q_{\ell}} . \]

\noindent (Note that if $q_{\ell}=1$ then $p_{\ell}$ might not be the integer closest to $q_{\ell} \alpha$.) Using the definition of $A$ and the fact that $\mathrm{sign} \, \varepsilon_{\ell} = (-1)^{\ell +1}$ from Proposition 5 (iv), we thus get that

\[ \left\| m \alpha \right\| = \left\| \frac{m p_{\ell}}{q_{\ell}} + \frac{m \varepsilon_{\ell}}{q_{\ell}} \right\| \ge \frac{1}{2} \left\| \frac{m p_{\ell}}{q_{\ell}} \right\| \]

\noindent for any $m \in A$. Hence

\[ \sum_{m \in A} \frac{1}{m^p \left\| m \alpha \right\|^p} = \sum_{a=1}^{\infty} \sum_{{{aq_{\ell} \le m < (a+1)q_{\ell}}\atop{m \in A}}} \frac{1}{m^p \left\| m \alpha \right\|^p} \le \sum_{a=1}^{\infty} \frac{1}{a^p q_{\ell}^p} \sum_{aq_{\ell} < m < (a+1)q_{\ell}} \frac{2^p}{\left\| \frac{m p_{\ell}}{q_{\ell}} \right\|^p} . \]

\noindent Since $p_{\ell}$ and $q_{\ell}$ are relatively prime, as $m$ runs in $a q_{\ell} < m < (a+1) q_{\ell}$, the numbers $m p_{\ell}$ fall into each nonzero residue class modulo $q_{\ell}$ exactly once. Thus the sum satisfies

\begin{equation}
\sum_{m \in A} \frac{1}{m^p \left\| m \alpha \right\|^p} \le \sum_{a=1}^{\infty} \frac{1}{a^p q_{\ell}^p} \sum_{j=1}^{q_{\ell}-1} \frac{2^p}{\left\| \frac{j}{q_{\ell}} \right\|^p} \le \sum_{a=1}^{\infty} \frac{2^{p+1}}{a^p} \zeta (p) = 2^{p+1} \zeta (p)^2 .
\end{equation}

Consider now the sum over $m \in B$. Taking the equation

\[ p_{\ell} q_{\ell-1} - q_{\ell}p_{\ell -1} = (-1)^{\ell} \]

\noindent from Proposition 5 (iii) modulo $q_{\ell}$, we learn that the multiplicative inverse of $p_{\ell}$ in the ring $\mathbb{Z}_{q_{\ell}}$ is $(-1)^{\ell} q_{\ell -1}$. This means that

\[ B = \left\{ a q_{\ell} + q_{\ell -1} : 1 \le a \le a_{\ell} -1 \right\} . \]

\noindent Indeed, the choice $a=a_{\ell}$ would result in $a_{\ell} q_{\ell} + q_{\ell -1} = q_{\ell +1}$ which is outside the interval $q_{\ell} \le m < q_{\ell +1}$. For any element $m = a q_{\ell} + q_{\ell -1} \in B$ we thus have

\[ \left\| m \alpha \right\| = \left\| \frac{(-1)^{\ell}}{q_{\ell}} + \frac{(a q_{\ell} + q_{\ell -1}) \varepsilon_{\ell}}{q_{\ell}} \right\| = \frac{1-q_{\ell -1} |\varepsilon_{\ell}|}{q_{\ell}} -a |\varepsilon_{\ell}| . \]

\noindent Rearranging the inequality

\[ \left| \varepsilon_{\ell} \right| = \left\| q_{\ell} \alpha \right\| \le \frac{1}{q_{\ell +1}} = \frac{1}{a_{\ell} q_{\ell} + q_{\ell -1}} \]

\noindent from Proposition 5 (i), we get

\[ a_{\ell} q_{\ell} | \varepsilon_{\ell} | \le 1 - q_{\ell -1} |\varepsilon_{\ell}| , \]

\noindent which in turn shows that for any $m=aq_{\ell} + q_{\ell -1} \in B$ we have

\[ \left\| m \alpha \right\| \ge (a_{\ell} -a) |\varepsilon_{\ell}| . \]

\noindent Therefore we have

\[ \sum_{m \in B} \frac{1}{m^p \left\| m \alpha \right\|^p} \le \sum_{a=1}^{a_{\ell}-1} \frac{1}{a^p q_{\ell}^p (a_{\ell}-a)^p |\varepsilon_{\ell}|^p} \le \frac{2^{p+1}}{a_{\ell}^p q_{\ell}^p |\varepsilon_{\ell}|^p} \zeta (p) . \]

\noindent Proposition 5 (i), (ii) imply

\begin{equation}
|\varepsilon_{\ell}| = \left\| q_{\ell} \alpha \right\| \ge \frac{1}{(a_{\ell}+1)q_{\ell} + q_{\ell -1}} \ge \frac{1}{3 a_{\ell} q_{\ell}} ,
\end{equation}

\noindent hence

\begin{equation}
\sum_{m \in B} \frac{1}{m^p \left\| m \alpha \right\|^p} \le 2 \cdot 6^p \zeta (p) .
\end{equation}

The sum over $m \in C$ will be the main term of (24). We have

\[ C = \left\{ a q_{\ell} : 1 \le a \le a_{\ell} \right\} , \]

\noindent since the choice $a=a_{\ell}+1$ would result in $(a_{\ell}+1)q_{\ell} > a_{\ell} q_{\ell} + q_{\ell -1} = q_{\ell +1}$. For any $m=a q_{\ell} \in C$ we have

\[ \left\| m \alpha \right\| = \left\| \frac{m \varepsilon_{\ell}}{q_{\ell}} \right\| = a |\varepsilon_{\ell}| , \]

\noindent therefore

\[ \sum_{m \in C} \frac{1}{m^p \left\| m \alpha \right\|^p} = \sum_{a=1}^{\infty} \frac{1}{a^{2p} q_{\ell}^p |\varepsilon_{\ell}|^p} - \sum_{a=a_{\ell}+1}^{\infty} \frac{1}{a^{2p} q_{\ell}^p |\varepsilon_{\ell}|^p} . \]

\noindent Using (26) again we obtain

\[ \left| \sum_{m \in C} \frac{1}{m^p \left\| m \alpha \right\|^p} - \frac{1}{q_{\ell}^p |\varepsilon_{\ell}|^p} \zeta (2p) \right| \le \int_{a_{\ell}}^{\infty} \frac{1}{x^{2p} q_{\ell}^p |\varepsilon_{\ell}|^p} \, \mathrm{d}x \le \frac{3^p}{2p-1} . \]

\noindent To estimate the main term, use Proposition 5 (i), (ii) to get

\[ \frac{1}{(a_{\ell}+2) q_{\ell}} \le |\varepsilon_{\ell}| = \left\| q_{\ell} \alpha \right\| \le \frac{1}{a_{\ell} q_{\ell}} , \]

\[ 0 \le \frac{1}{q_{\ell}^p |\varepsilon_{\ell}|^p} - a_{\ell}^p \le (a_{\ell}+2)^p - a_{\ell}^p \le 2 p 3^{p-1} a_{\ell}^{p-1} . \]

\noindent Hence we have

\begin{equation}
\left| \sum_{m \in C} \frac{1}{m^p \left\| m \alpha \right\|^p} - \zeta (2p) a_{\ell}^p \right| \le \frac{3^p}{2p-1} + 2 p 3^{p-1} \zeta (2p) a_{\ell}^{p-1} .
\end{equation}

\noindent (25), (27) and (28), together with the trivial estimate

\[ 2^{p+1} \zeta (p)^2 + 2 \cdot 6^p \zeta (p) + \frac{3^p}{2p-1} + 2 p 3^{p-1} \zeta (2p) < 6^p \frac{4p^2}{(p-1)^2} \]

\noindent show that

\begin{equation}
\left| \sum_{q_{\ell} \le m < q_{\ell}+1} \frac{1}{m^p \left\| m \alpha \right\|^p} - \zeta (2p) a_{\ell}^p \right| < 6^p \frac{4 p^2}{(p-1)^2} \cdot a_{\ell}^{p-1}
\end{equation}

\noindent holds whenever $\ell \ge 2$ and $q_{\ell} \ge 2$.

We now claim that (29) in fact holds for every $\ell \ge 1$. We have to distinguish between two cases. If $a_1 >1$, then $1=q_1 < q_2 < \dots$, thus we need to prove (29) for $\ell =1$ only. If, on the other hand, $a_1=1$, then $1=q_1=q_2 < q_3 < \dots$, thus we have to prove (29) for $\ell = 1,2$.

Suppose that $a_1 >1$. Recalling the algorithm for finding the continued fraction representation of $\alpha$, we have

\[ \frac{1}{a_1 +1} \le \alpha - a_0 \le \frac{1}{a_1} . \]

\noindent Hence for any $1 \le m \le \frac{a_1}{2}$ we have

\[ \frac{m}{a_1 +1} \le \left\| m \alpha \right\| = m \alpha - m a_0 \le \frac{m}{a_1} , \]

\[ 0 \le \frac{1}{m^p \left\| m \alpha \right\|^p} - \frac{a_1^p}{m^{2p}} \le \frac{(a_1+1)^p - a_1^p}{m^{2p}} \le \frac{p 2^{p-1}}{m^{2p}} a_1^{p-1} .  \]

\noindent On the other hand, if $\frac{a_1+1}{2} \le m \le a_1 -1$, then

\[ \frac{a_1-m}{a_1} \le \left\| m \alpha \right\| = \left| m \alpha - m a_0 -1 \right| \le \frac{a_1+1-m}{a_1+1}. \]

\noindent Therefore

\[ \left| \sum_{q_1 \le m < q_2} \frac{1}{m^p \left\| m \alpha \right\|^p} - \zeta (2p) a_1^p \right| \le \]

\[ \sum_{1 \le m \le \frac{a_1}{2}} \frac{p 2^{p-1}}{m^{2p}} a_1^{p-1} + \sum_{m \ge \frac{a_1+1}{2}} \frac{a_1^p}{m^{2p}} + \sum_{\frac{a_1+1}{2} \le m \le a_1-1} \frac{a_1^p}{m^p (a_1-m)^p} \le 6^p \frac{4 p^2}{(p-1)^2} a_1^{p-1} . \]

Finally, suppose that $a_1=1$. Then (29) is clearly true for $\ell=1$. We have $q_1=q_2=1$ and $q_3 = a_2+1$, thus we need to consider

\[ \sum_{q_2 \le m < q_3} \frac{1}{m^p \left\| m \alpha \right\|^p} = \sum_{m=1}^{a_2} \frac{1}{m^p \left\| m \alpha \right\|^p} . \]

\noindent Recalling the algorithm for finding the continued fraction representation of $\alpha$, we have

\[ \frac{1}{a_2+1} \le \frac{1}{\alpha - a_0} -1 \le \frac{1}{a_2} , \]

\[ \frac{-1}{a_2+1} \le \alpha - a_0 -1 \le \frac{-1}{a_2+2} . \]

\noindent Hence for any $1 \le m \le \frac{a_2+1}{2}$ we have

\[ \frac{m}{a_2 +2} \le \left\| m \alpha \right\| = \left| m \alpha - m a_0 -m \right| \le \frac{m}{a_2 +1} , \]

\[ 0 \le \frac{1}{m^p \left\| m \alpha \right\|^p} - \frac{a_2^p}{m^{2p}} \le \frac{(a_2+2)^p - a_2^p}{m^{2p}} \le \frac{2p3^{p-1} a_2^{p-1}}{m^{2p}} . \]

\noindent On the other hand, if $\frac{a_2 +2}{2} \le m \le a_2$, then

\[ \frac{a_2+1-m}{a_2+1} \le \left\| m \alpha \right\| = m \alpha - m a_0 -m+1 \le \frac{a_2+2-m}{a_2+2} . \]

\noindent Therefore

\[ \left| \sum_{q_2 \le m < q_3} \frac{1}{m^p \left\| m \alpha \right\|^p} - \zeta (2p) a_2^p \right| \le \]

\[ \sum_{1 \le m \le \frac{a_2+1}{2}} \frac{2p 3^{p-1} a_2^{p-1}}{m^{2p}} + \sum_{m \ge \frac{a_2+2}{2}} \frac{a_2^p}{m^{2p}} + \sum_{\frac{a_2+2}{2} \le m \le a_2} \frac{(a_2+1)^p}{m^p (a_2 +1-m)^p} \le 6^p \frac{4p^2}{(p-1)^2} a_2^{p-1}. \]

\noindent This finishes the proof of the fact that (29) holds for every integer $\ell \ge 1$.

Summing (29) we obtain

\[ \left| \sum_{0 < m < q_{\ell}} \frac{1}{m^p \left\| m \alpha \right\|^p} - \sum_{0 < k < \ell} \zeta (2p) a_k^p \right| \le 6^p \frac{4p^2}{(p-1)^2} \sum_{0<k<\ell} a_k^{p-1} . \]

\noindent Using the general inequality $|A-B| \le \frac{\left| A^p - B^p \right|}{B^{p-1}}$, which holds for any positive reals $A,B$, we get

\[ \left| \left( \sum_{0 < m < q_{\ell}} \frac{1}{m^p \left\| m \alpha \right\|^p} \right)^{\frac{1}{p}} - \left( \sum_{0 < k < \ell} \zeta (2p) a_k^p \right)^{\frac{1}{p}} \right| \le 6^p \frac{4p^2}{(p-1)^2} \cdot \frac{\sum_{0<k<\ell} a_k^{p-1}}{\left( \sum_{0<k<\ell} \zeta (2p) a_k^p \right)^{\frac{p-1}{p}}} . \]

\noindent Finally, the inequality

\[ \left( \frac{1}{\ell -1} \sum_{0<k< \ell} a_k^{p-1} \right)^{\frac{1}{p-1}} \le \left( \frac{1}{\ell -1} \sum_{0<k < \ell} a_k^p \right)^{\frac{1}{p}} \]

\noindent yields

\[ 6^p \frac{4p^2}{(p-1)^2} \frac{\sum_{0<k<\ell} a_k^{p-1}}{\left( \sum_{0<k<\ell} \zeta (2p) a_k^p \right)^{\frac{p-1}{p}}} \le 6^p \frac{4p^2}{(p-1)^2} (\ell -1)^{\frac{1}{p}} . \]

\begin{flushright}
$\square$
\end{flushright}

\vspace{10mm}

\noindent\textbf{References}

\begin{enumerate}
\item[{[1]}] P. Barkan. \textit{Sur les sommes de Dedekind et les fractions continues finies.} C. R. Acad. Sci. Paris S\'er. A-B 284 (1977), no. 16, A923-A926.

\item[{[2]}] J. Beck. \textit{Probabilistic Diophantine approximation. Randomness in lattice point counting.} Springer Monographs in Mathematics. Springer, Cham, 2014. xvi+487 pp. ISBN: 978-3-319-10740-0.

\item[{[3]}] D. Bilyk, V. N. Temlyakov, R. Yu. \textit{Fibonacci sets and symmetrization in discrepancy theory.} J. Complexity 28 (2012), no. 1, 18-36.

\item[{[4]}] J. W. S. Cassels. \textit{An introduction to Diophantine approximation.} Cambridge Tracts in Mathematics and Mathematical Physics, no. 45. Cambridge University Press, New York, 1957. x+166 pp.

\item[{[5]}] H. Davenport. \textit{Note on irregularities of distribution.} Mathematika 3 (1956), 131-135.

\item [{[6]}] R. R. Hall, J. C. Wilson. \textit{On reciprocity formulae for inhomogeneous and homogeneous Dedekind sums.} Math. Proc. Cambridge Philos. Soc. 114 (1993), no. 1, 9-24.

\item[{[7]}] R. Kritzinger, L. M. Kritzinger. \textit{$L_2$ discrepancy of symmetrized generalized Hammersley point sets in base $b$.} J. Number Theory 166 (2016), 250-275.

\item[{[8]}] K. Roth. \textit{On irregularities of distribution.} Mathematika 1 (1954), 73-79.

\end{enumerate}

\end{document}